\documentclass[11pt]{amsart}

\usepackage{amsmath}

\usepackage{amssymb}

\newtheorem{thm}{Theorem}[section]
\newtheorem{prop}[thm]{Proposition}
\newtheorem{lemma}[thm]{Lemma}
\newtheorem{cor}[thm]{Corollary}

\theoremstyle{definition}
\newtheorem{defn}[thm]{Definition}
\newtheorem{ex}[thm]{Example}

\theoremstyle{remark}
\newtheorem{remark}[thm]{Remark}

\numberwithin{equation}{section}

\def\R{\mathbb{R}}
\def\RP{\mathbb{RP}}

\def\Z{\mathbb{Z}}

\def\Q{\mathcal{Q}}
\def\D{\mathcal{D}}
\def\P{\mathcal{P}}

\def\dim{\mathrm{dim\,}}
\def\codim{\mathrm{codim\,}}

\def\span{\mathrm{span\,}}

\begin{document}

\title[Conformal deformation]{Conformal deformation of spacelike\\
surfaces in Minkowski space}

\author{Emilio Musso}
\address{(E. Musso) Dipartimento di Matematica Pura ed Applicata,
Universit\`a degli Studi dell'Aquila, Via Vetoio, I-67010
Coppito (L'Aquila), Italy} \email{musso@univaq.it}

\author{Lorenzo Nicolodi}
\address{(L. Nicolodi) Di\-par\-ti\-men\-to di Ma\-te\-ma\-ti\-ca,
Uni\-ver\-si\-t\`a degli Studi di Parma, Viale G. P. Usberti 53/A -
Campus universitario, I-43100 Parma, Italy} \email{lorenzo.nicolodi@unipr.it}

\thanks{Authors partially supported by MIUR projects:
\textit{Metriche riemanniane e variet\`a differenziali} (E.M.);
\textit{Propriet\`a geometriche delle variet\`a reali e complesse} (L.N.);
and by the GNSAGA of INDAM}





\begin{abstract}
We address the problem of second order conformal deformation of
spacelike surfaces in compactified Minkowski 4-space. We explain the
construction of the exterior differential system of conformal
deformations and discuss its general and singular solutions. In
particular, we show that isothermic surfaces are singular solutions
of the system, which implies that a generic second order deformable
surface is not isothermic.
This differs from the situation in 3-dimensional conformal geometry, where 
isothermic surfaces coincide with
deformable surfaces.
\end{abstract}

\maketitle

\section{Introduction}\label{s:intro}

The surfaces in the conformal 3-sphere which admit second order
deformations with respect to the group of conformal transformations
coincide with isothermic surfaces \cite{Ca3}, \cite{Mtrieste}. This
is no longer true in higher dimensional conformal spaces. Actually,
isothermic surfaces are deformable to second order \cite{HJlibro},
but generically a deformable surface is not isothermic. This result
was originally stated without proof by E. Cartan in his address at
the 1920 International Congress of Mathematicians \cite{Ca1}. More
precisely, as an illustration of the general deformation theory of
submanifolds in homogeneous spaces, Cartan indicated that isothermic
surfaces in conformal 4-space are singular solutions of the exterior
differential system (EDS) which defines deformable surfaces.
The work of Cartan on the deformation of submanifolds
in homogeneous spaces and the related questions of contact and rigidity
were taken up and further developed by P. Griffiths and
G. Jensen  \cite{Gr}, \cite{J}.
In particular, that the problems of $k$th order deformation are
equivalent to solving certain EDSs on appropriate spaces of frames was
established in \cite{J}.
Still, for each concrete geometric situation there is a specific
problem to solve.

\vskip0.1cm

In this paper, we investigate the problem of conformal deformation
for the case of spacelike surfaces in compactified Minkowski
4-space.\footnote{In the case of spacelike surfaces, the discussion
is independent of the signature of the target space one considers.}
We give a detailed description of the Pfaffian differential system
(PDS) of conformal deformations and then provide the tools to
discuss its general and singular solutions within the theory of
EDSs.\footnote{We recall that an integral manifold $f : N \to M$ of
an EDS $\mathcal{I}$ on $M$ is a \textit{general solution} if the
integral element $df(T_p N)$ is ordinary, for every $p\in N$; it is
a \textit{singular solution} if the integral element $df(T_p N)$
fails to be ordinary, for every $p\in N$. So, singular solutions are
not given by the Cartan--K\"ahler theorem. They involve additional
equations (see \cite{BCGGG} for more detail).} We specify the
appropriate space which supports the differential system of a
conformal deformation and find the equations of the integral
elements of the system. We show that the differential system of
deformations is in involution and that its general solutions depend
on one arbitrary function in two variables.
We then determine the equations of the variety defining the singular
solutions of the system and show that isothermic surfaces are indeed
singular solutions. In particular, isothermic surfaces depend on six
arbitrary functions in one variable, which implies that a generic
second order deformable surface is not isothermic.

\vskip0.1cm

The paper is organized as follows. In Section \ref{s:pre}, we
present the background material and set up the basic constructions.
In Section \ref{s:def}, we discuss the questions of conformal
deformation and rigidity for spacelike surfaces in compactified
Minkowski space and describe isothermic surfaces as examples of
deformable surfaces. In Section \ref{s:eds}, we study the
involutiveness of the PDS of a conformal deformation and discuss its
general and singular solutions.

\vskip0.1cm

We use \cite{BCGGG} as basic reference for the theory of EDSs.
For a general account on submanifold theory in conformal differential
geometry, we refer to \cite{HJlibro}. The summation convention on
repeated indices is used throughout the paper.

\vskip0.1cm

\noindent \textit{Acknowledgments.} The authors would like to thank the
referee, whose comments and suggestions greatly
contributed to improve the first version of the paper.

\section{Preliminaries and basic constructions}\label{s:pre}

\subsection{The conformal completion of Minkowski 4-space}

Let $\R^{4,2}$ denote $\R^6$ with the symmetric bilinear form
\begin{equation}\label{scalar-product}
 \langle x,y\rangle =
   - (x^0y^5 + x^5y^0) + x^1y^1 + x^2y^2 + x^3y^3 - x^4y^4 = g_{IJ}x^Iy^J
    \end{equation}
of signature $(4,2)$, where $x^0,\dots,x^5$ are the coordinates with
respect to the standard basis $e_0,\dots,e_5$ of $\R^{6}$. The \textit{Lie
quadric} is the hypersurface
\[
 \mathcal{Q} = \{[k]\in \RP^5 \, | \, \langle k, k\rangle=0\}.
  \]
A pair of points $[k_1], [k_2]$ in $\mathcal{Q}$ satisfying $\langle
k_1, k_2\rangle=0$ defines a line $[k_1,k_2]$ in $\mathcal{Q}$. The
set of all lines in $\mathcal{Q}$ forms a smooth manifold of
dimension $5$, which we denote by $\Lambda$. Let $\mathrm{O}(4,2)$
denote the pseudo-orthogonal group of (\ref{scalar-product}). The
standard action of $\mathrm{O}(4,2)$ on $\RP^5$ maps the quadric
$\mathcal{Q}$ into itself and induces an action on $\Lambda$ which
is transitive.


 The quadric $\mathcal{Q}$ is diffeomorphic to $(S^1\times
S^3)/{\Z_2}$ and inherits a locally conformally flat metric of
signature $(3,1)$ from the flat metric on $\R^{4,2}$ corresponding
to \eqref{scalar-product}. This implies that
$\mathrm{O}(4,2)/{\Z_2}$ acts on $\mathcal{Q}$ as a group of
conformal transformations. In fact, $\mathrm{O}(4,2)/{\Z_2}$
coincides with the conformal group of $\mathcal{Q}$. The Lie quadric
$\mathcal{Q}$ can be regarded as the conformal compactification of
flat Minkowski spacetime $\R^{3,1}$ (see \cite{K1949},
\cite{GSkepler}). The conformal embedding of $\R^{3,1}$ with flat
Lorentz metric $(\cdot,\cdot)$ is given by
\[
 \R^{3,1} \ni v
  \mapsto
   \biggl[\Bigl(1,v,\frac{1}{2}(v,v)\Bigr)^T \biggr] \in \mathcal{Q}.
     \]

\begin{remark}[Lie sphere geometry]\label{r:lie-geom}
The points in the Lie quadric $\mathcal{Q}$ are in bijective
correspondence with the set of all oriented spheres and point
spheres in the unit sphere $S^3 \subset\R^4$. A line $[k_1,k_2]$ in
$\mathcal{Q}$ corresponds to a family of spheres in oriented
contact. This family of oriented spheres contains a unique point
sphere, which is the common point of contact, and determines the
common unit normal vector at this point. The set $\Lambda$ can then
be identified with $T_1S^3=\{(u,v)\in S^3\times S^3 \subset
\R^4\times \R^4 \, |  \, u\cdot v=0\}$. A Lie sphere transformation
is a projective transformation induced by a transformation in the
group $\mathrm{O}(4,2)/{\Z_2}$. In terms of $S^3$, a Lie sphere
transformation in a map on the space of oriented spheres which
preserves oriented contact. The group $\mathrm{O}(4,2)/{\Z_2}$ acts
on $\Lambda$, and hence on $T_1S^3$, as a group of contact
transformations. See \cite{Ce} for more detail.
\end{remark}


We now introduce moving frames to study surface theory in
$\mathcal{Q}$. Let $G$ be the identity component of
$\mathrm{O}(4,2)$ and $\mathfrak{g} = \{B\in \mathfrak{gl}(6,\R) \,
|\, {B^T}g+gB=0 \}$ its Lie algebra, where $g=(g_{IJ})$. By a frame
is meant a basis $A_0,\dots,A_5$ of $\R^{4,2}$ such that
$(A_0,\dots,A_5)\in G$. Up to the choice of a reference frame, the
manifold of frames identifies with $G$. For $A\in G$, let $A_J =
Ae_J$ be the column vectors of $A$ and regard the $A_J$ as
$\R^{4,2}$ valued functions on $G$. Since the $A_J$ form a basis of
$\R^{4,2}$, there exist unique left invariant 1-forms $\omega^I_J$
($I,J = 0,1,\dots,5$) such that
\begin{equation}\label{dA}
  dA_J=\omega^I_JA_I \quad (J=0,\dots,5).
   \end{equation}
Differentiating \eqref{dA} yields the Cartan structure equations
\[
   d\omega^I_J = - \omega^I_K \wedge \omega^K_J.
   \]
Differentiating $\langle A_I,A_J\rangle= g_{IJ}$ gives the symmetry
equations
\[
   \omega^K_Ig_{KJ} + \omega^K_Jg_{KI} = 0.
   \]
The 1-forms $\omega^I_J$ are the components of the {Maurer--Cartan
form} $\omega=A^{-1}dA$ of $G$, which accordingly takes the form
\begin{equation}\label{MCform}
\omega=\left(
\begin{array}{cccccc}
\omega^0_0&\omega^0_1&\omega^0_2&\omega^0_3&\omega^0_4&0\\
\omega^1_0&0&-\omega^2_1&-\omega^3_1&\omega^4_1&\omega^0_1\\
\omega^2_0&\omega^2_1&0&-\omega^3_2&\omega^4_2&\omega^0_2\\
\omega^3_0&\omega^3_1&\omega^3_2&0&\omega^4_3&\omega^0_3\\
\omega^4_0&\omega^4_1&\omega^4_2&\omega^4_3&0&-\omega^0_4\\
0&\omega^1_0&\omega^2_0&\omega^3_0&-\omega^4_0&-\omega^0_0
\end{array}\right).
\end{equation}

The standard action of $G$ on $\RP^5$ restricts to a transitive
action on the quadric $\mathcal{Q}$. This action defines a principal
$K_{\mathcal{Q}}$-bundle
\[
 \pi_\mathcal{Q} : G \to \mathcal{Q}\cong G/K_{\mathcal{Q}},
  \quad A \mapsto A[e_0]=[A_0],
   \]
where $K_{\mathcal{Q}}$ is the isotropy subgroup at $[e_0]$. It is
easy to compute that
\[
K_{\mathcal{Q}} = \left\{\left(
 \begin{array}{ccc}
  r&y^TB&\displaystyle\frac{1}{2r}(y,y)\\
  0&B&\displaystyle {y}/{r}\\
  0&0&\displaystyle {1}/{r}\\
  \end{array}\right)
\,: \,
B \in \mathrm{SO}(3,1),\,
  y\in \R^4,\, r>0
 \right\}.
     \]
From this it follows that the forms $\{\omega^1_0, \omega^2_0,
\omega^3_0, \omega^4_0\}$ span the semibasic forms of the projection
$\pi_\mathcal{Q}$.\footnote{We recall that a differential form
$\varphi$ on the total space of a fiber bundle $\pi : P \to B$ is
said to be \textit{semibasic} if its contraction with any vector
field tangent to the fibers of $\pi$ vanishes, or equivalently, if
its value at each point $p\in P$ is the pullback via $\pi^\ast_p$ of
some form at $\pi(p)\in B$. Some authors call such a form
\textit{horizontal}. A stronger condition is that $\varphi$ be
\textit{basic}, meaning that it is locally the pullback via
$\pi^\ast$ of a form on the base $B$.} The conformal structure on
$\mathcal{Q}$ is determined by the quadratic form
\[
  (\omega^1_0)^2 + (\omega^2_0)^2 + (\omega^3_0)^2  -
   (\omega^4_0)^2.
    \]

\vskip0.2cm

Next, we will introduce two additional
homogeneous spaces which will play a role in the discussion of 
the conformal deformation problem.

\subsection{The Grassmannian of parabolic 3-planes}

Let $\mathcal{P}\subset G_3(\R^{4,2})$ be the set of degenerate
3-planes $V\subset\R^{4,2}$ of signature $(0,+, +)$, i.e.,

\begin{itemize}
\item $\dim V\cap V^\perp = 1$;

\item $\langle \cdot,\cdot\rangle$ restricted to $V/V\cap V^\perp$ is positive definite.
\end{itemize}

\noindent The natural action of $G$ on $\mathcal{P}$ is transitive and the
map
\[
  \pi_\mathcal{P} : G \to \mathcal{P},\quad A \mapsto A\cdot [e_0\wedge e_1\wedge e_2]= [A_0\wedge A_1\wedge A_2]
   \]
makes $G$ into a principal fiber bundle over $\P$ with fiber
\[
 H_\mathcal{P} = \left\{ A \in G \,: \; A\cdot [e_0\wedge e_1\wedge e_2] =
  [e_0\wedge e_1\wedge e_2]\right\},
   \]
where $H_\mathcal{P} \subset G$ is the isotropy subgroup at $[e_0\wedge e_1\wedge e_2]$.
A direct computation shows that $H_\mathcal{P}$ consists of matrices of the form
\begin{equation}\label{isotropy1}
  A(r,x,y,a,b) =
 \left(
 \begin{array}{cccc}
  r&x^Ta&y^T\mathbb{I}_{1,1}b&\displaystyle\frac{1}{2r}(x^Tx + y^T\mathbb{I}_{1,1}y)\\
  0&a&0&\displaystyle {x}/{r}\\
  0&0&b&\displaystyle {y}/{r}\\
  0&0&0&\displaystyle {1}/{r}\\
  \end{array}\right),
     \end{equation}
where $a \in \mathrm{SO}(2)$, $b \in \mathrm{SO}(1,1)$, $x,y\in \R^2$, $r>0$, and
\[
 \mathbb{I}_{1,1}= \left(
  \begin{array}{cc}
   1&0\\
    0&-1
    \end{array}\right).
    \]
This implies that $\mathcal{P}$ is an 8-dimensional
homogeneous space of $G$ and that
$$
  \{\omega^1_0,\, \omega^2_0,\, \omega^3_0,\, \omega^4_0,\, \omega^3_1,\,\omega^3_2,\,
   \omega^4_1,\, \omega^4_2\}
    $$
is a basis for the space of semibasic 1-forms of the projection $\pi_\mathcal{P}$.

\subsection{The configuration space and some relevant PDSs}
\label{ss:conf-space}

Let $\mathcal{D}$ be the submanifold of
$\mathcal{P}\times\mathcal{P}\times G$ defined by
\[
  \mathcal{D} := \left\{ (V_1,V_2,A) \in \mathcal{P}\times\mathcal{P}\times G\,:\;
    A\cdot V_1 = V_2  \right\}.
     \]
We call $\mathcal{D}$ the \textit{configuration space of
deformations}. The Lie group $G\times G$ acts on the left
on $\mathcal{D}$ by
\[
  (A,B)\cdot (V_1,V_2,F) := \left(A\cdot V_1, B\cdot V_2, B F A^{-1}\right),
  \]
for each $(A,B) \in G\times G$, and $(V_1,V_2,F) \in \mathcal{D}$.
This action is transitive and the isotropy subgroup of $G\times G$ at
\[
  \left([e_0\wedge e_1\wedge e_2],[e_0\wedge e_1\wedge e_2], {e}_{G} \right) \in \mathcal{D}
    \] 
(where ${e}_G$ is the identity element of $G$) is the closed subgroup
\[
   {(G\times G)}_\mathcal{D} = \left\{(A,B) \in {G\times G} \, : \; A=B \in H_\mathcal{P}
     \right\}.
      \]
Thus $\mathcal{D}$ is a 23-dimensional homogeneous space of
${G\times G}$ and the natural projection $\pi_\mathcal{D} : G\times G
\to \mathcal{D}$ is given by
\[
   (A,B) \mapsto \left([A_0\wedge A_1\wedge A_2],[B_0\wedge B_1\wedge B_2],BA^{-1}\right).
    \]

Next, let $(\omega, \Omega)$ denote the Maurer--Cartan form of ${G\times G}$.
The set of left-invariant forms
\begin{equation}\label{basicforms}
  \omega^0_0 - \Omega^0_0, \, \omega^0_I - \Omega^0_I,\, \omega^2_1 - \Omega^2_1,\,
  \omega^4_3 - \Omega^4_3,\, \omega^a_i,\, \Omega^a_i,\, \omega^I_0,\, \Omega^I_0,
     \end{equation}
where $I = 1,2,3,4$; $i= 1,2$; $a=3,4$, is a basis for the space of
semibasic forms of the fibration $\pi_\mathcal{D}$.
Let
\[
  \alpha^0_0 - \beta^0_0, \, \alpha^0_I - \beta^0_I,\, \alpha^2_1 - \beta^2_1,\,
  \alpha^4_3 - \beta^4_3,\, \alpha^a_i,\, \beta^a_i,\, \alpha^I_0,\, \beta^I_0
     \]
be the 1-forms on $\mathcal{D}$ obtained
by pulling back the forms \eqref{basicforms}
via a local section of $\pi_\mathcal{D}$.

Out of these forms, we can construct three
invariant Pfaffian systems $I_1$, $I_2$, $I_3\subset \Gamma (T^\ast\mathcal{D})$ given by
\begin{gather*}
 I_1 =\span\{\alpha^i_0 - \beta^i_0,\, \alpha^a_0,\,\beta^a_0 \},\\
   I_2 = \span\{\alpha^i_0 - \beta^i_0,\, \alpha^a_0, \,\beta^a_0,\,
         \alpha^0_0 - \beta^0_0,\, \alpha^2_1 - \beta^2_1, \,\alpha^a_i - \beta^a_i \},\\
      I_3 = \span \{\alpha^i_0 - \beta^i_0,\, \alpha^a_0,\, \beta^a_0, \,
             \alpha^0_0 - \beta^0_0,\, \alpha^2_1 - \beta^2_1,\, \alpha^a_i - \beta^a_i,\,
             \alpha^4_3 - \beta^4_3,\, \alpha^0_a - \beta^0_a \},
  \end{gather*}
$i= 1,2$; $a=3,4$.
In Section \ref{s:def}, we will explain the geometric meaning of
the Pfaffian differential systems defined by the differential ideals
$\mathfrak{I}_1$, $\mathfrak{I}_2$, and $\mathfrak{I}_3$
generated by $I_1$, $I_2$, and $I_3$, respectively.

\subsection{Spacelike surfaces in the Lie quadric}

Let $X$ be a 2-dimensional oriented manifold and let
$f : X \to \Q \subset \RP^5$ be a smooth spacelike conformal immersion.

\begin{defn}
A \textit{zeroth order frame field along} $f$ is a smooth map $A :
U \subset X \to G$, defined on an open subset $U\subset X$, such that
$f = [A_0] = \pi_\mathcal{Q}\circ A$.
\end{defn}

For any such a frame we put
$\theta= A^\ast \omega$. The totality of zeroth order frames along
$f$ is the bundle $\pi_0 : P_0(f) \to X$, where
\[
  P_0(f) = \left\{(q,A) \in X\times G \, : \, [A_0] = f(q) \right\}.
   \]

\begin{defn}
A \textit{first order frame field along} $f$ is a zeroth order frame field
$A : U\subset X\to G$ such that
\[
 \theta^3_0 = \theta^4_0 = 0.
  \]
The totality of first order frames gives rise to a subbundle
\[
  \pi_1 : P_1(f) \to X
   \]
of $P_0(f)$, referred to as the \textit{first order frame bundle of} $f$.
The structure group of $P_1(f)$ consists of matrices of the form \eqref{isotropy1}.
\end{defn}

The map
$$
  P_1(f) \ni (q,A) \mapsto [A_0\wedge A_1 \wedge A_2] \in \P
  $$
is constant along the fibers of $\pi_1 : P_1(f) \to X$, and therefore induces a
well-defined map
$$
 \tau_f : X \to \P
   $$
and a corresponding rank 3 vector bundle
$$
  \tau(X) = \left\{(q,W)\in X\times \R^{4,2} \, : \, W \in \tau_f(q)\right\} \to X.
  $$
The tautological line bundle
$$
  K_X = \left\{(q,W)\in X\times \R^{4,2} \, : \, W \in f(q)\right\} \to X
  $$
is a line subbundle of $\tau(X)$ such that $K_X = \tau(X)\cap
{\tau(X)}^\perp$. Thus the quotient bundle
 $$
  \mathcal{T}(X) = \tau(X)/K_X \to X
   $$
inherits from $\R^{4,2}$ a Riemannian metric, say $g^\tau$. Note
that the tensor product $\mathcal{T}(X) \otimes K_X^\ast$ can be
canonically identified with the tangent bundle $T(X)$.

Similarly, the map
$$
  P_1(f) \ni (q,A) \mapsto [A_0\wedge A_3 \wedge A_4]
  $$
is constant along the fibers of $\pi_1 : P_1(f) \to X$, and induces a
well-defined map
$$
  \nu_f : X \to G_3(\R^{4,2}).
  $$
We denote by $\nu(X) \to X$ the corresponding rank 3 vector bundle
$$
  \nu(X) = \left\{(q,W)\in X\times \R^{4,2} \, : \, W \in \nu_f(q)\right\} \to X.
  $$
Again, $K_X = \nu(X)\cap {\nu(X)}^\perp$ and the quotient bundle
$$
  \mathcal{N}(X) = \nu(X)/K_X \to X
   $$
inherits from $\R^{4,2}$ a pseudo-Riemannian metric $g^\nu$ of
signature $(1,1)$. We call $\mathcal{N}(X)$ the \textit{conformal
normal bundle} of the spacelike immersion $f : X \to \Q$. The
conformal normal bundle is equipped with a metric covariant
derivative $D^\nu$ defined by
$$
  D^\nu (x\tilde{A}_3 + y \tilde{A}_4) = (dx + y \theta^4_3)\tilde{A}_3 +
   (dy + x\theta^4_3)\tilde{A}_4,
    $$
where $A : U \to G$ is a first order frame and
$(\tilde{A}_3,\tilde{A}_4)$ denotes the induced local trivialization
of the conformal normal bundle. Note that $g^\tau$ induces a Riemannian
metric on the symmetric tensor product
$S^2\mathcal{T}^\ast(X)$.
Let $A : U \to G$ be a fixed first order frame along $f$ and denote
by $\tilde{A}_0$, $(\tilde{A}_1,\tilde{A}_2)$ and
$(\tilde{A}_3,\tilde{A}_4)$ the corresponding trivializations of the
bundles $K_X$, $\mathcal{T}(X)$ and $\mathcal{N}(X)$, respectively.
Differentiating $\theta^3_0 =\theta^4_0 = 0$ and applying Cartan's
Lemma yield
$$
 \theta^a_i = h^a_{ij} \theta^j_0, \quad  h^a_{ij} = h^a_{ji}
   \quad (a=3,4; i,j = 1,2)
     $$
for smooth functions $h^a_{ij} : U \to \R$.

The \textit{conformal second fundamental form} of $f$ is the trace-free quadratic
form $\mathcal{A}$, taking values in $\mathcal{N}(X)\otimes K_X$,
given locally by
\[
 \mathcal{A} = \biggl[h^a_{ij} - \Bigl(\sum_{l=1,2}h^a_{ll}\Bigr)\delta_{ij}\biggr]
  \theta^i_0  \theta^j_0 \otimes A_a \otimes A_0.
   \]
The form $\mathcal{A}$
is independent of the choice of first order frames and is a conformal
invariant of the immersion $f$.

\begin{defn}
A \textit{second order frame field along} $f$ is a first order frame field
$A : U \subset X \to G$ such that
\[
 \sum_{l=1,2} h^a_{ll}=0 \quad (a=3,4).
  \]
\end{defn}
Locally there exist second order frames. The totality of second
order frame fields gives rise to a subbundle
$\pi_2 : P_2(f) \to X$
of $P_1(f)$, referred to as the \textit{second order frame bundle of} $f$, whose
structure group is the subgroup $G_2 \subset G$ given by
\[
  G_2 = \left\{
\left(
 \begin{array}{cccc}
  r&x^Ta&0&\displaystyle\frac{x^Tx}{2r}\\
  0&a&0&\displaystyle {x}/{r}\\
  0&0&b&0\\
  0&0&0&\displaystyle {1}/{r}\\
  \end{array}\right) \, :
\begin{array}{ll}
  a \in \mathrm{SO}(2),&b \in \mathrm{SO}(1,1),\\
  x\in \R^2,&r>0
  \end{array}
  \right\}.
     \]

\begin{remark}
The frame bundles associated with an immersed surface in $\Q$
considered above are the analogs of the bundles considered by Bryant
for an immersed surface in the conformal 3-sphere \cite{BrDG}.
\end{remark}

\begin{remark}
Observe that the mappings
\[
  [A_3 -A_4], [A_3 + A_4] : X \to \Q
  \]
are independent of the choice of the second order frame $A$. Further,
\begin{eqnarray*}
 &&F_1 = [A_0 \wedge (A_3 - A_4)] : X \to \Lambda,\\
 &&F_2 = [A_0 \wedge (A_3 + A_4)] : X \to \Lambda
   \end{eqnarray*}
are two Legendrian immersions with respect to the canonical contact structure
of $\Lambda$ (see Remark \ref{r:lie-geom}). If
$p : \Lambda \to \R^3 \cup \{\infty\}$
denotes the projection of $\Lambda$ onto the 3-sphere, the mappings
\begin{eqnarray*}
  &&\phi_1 = p \circ F_1 : X \to \R^3 \cup \{\infty\}. \\
  &&\phi_2 = p \circ F_2 : X \to \R^3 \cup \{\infty\}
   \end{eqnarray*}
are the two envelopes of the congruence of spheres represented by
the spacelike immersion $f$ (see Remark \ref{r:lie-geom}). Note that
$(A_3 - A_4)$ and $(A_3 + A_4)$ generate the isotropic line
subbundle of the conformal normal bundle $\mathcal{N}(X)$.
\end{remark}

\section{Contact and deformation}\label{s:def}

We start by recalling the notion of deformation (see \cite{Gr}, \cite{J}).

\begin{defn}
Let $G/H$ be a homogeneous space and let $f,\hat{f} : X\to G/H$ be
smooth maps. Two such maps $f$ and $\hat{f}$ are $k$th \textit{order
$G$-deformations of each other} if there exists a smooth map ${D} :
X \to G$ such that, for each point $p\in X$, the maps $\hat{f}$ and
${D}(p)f$ have analytic contact of order $k$ at $p$; that is, if
they have the same $k$th order jets at $p$. The map ${D}$ is called
the \textit{infinitesimal displacement} of the deformation. When
$D(p)$ does not depend on $p\in X$, the deformation is called
\textit{trivial}, and then $\hat{f} = Df$ is $G$-congruent to $f$. A
given map $f:X\to G/H$ is \textit{rigid} to $k$th order deformations
if there are no nontrivial $k$th order deformations of it; it is
\textit{deformable of order} $k$ if it admits a nontrivial $k$th
order deformation.

\end{defn}

\subsection{Analytic Contact}

Let $X$ be a 2-dimensional manifold and $f,\,\hat{f} : X \to \Q$
smooth maps. To express the condition of $k$th order analytic
contact for $f$ and $\hat{f}$, we need to introduce some notation.
Let $\{U; x^1,x^2\}$ be a local coordinate system, where $U\subset
X$ is an open set. Let $S^h(U)$ be the space of symmetric $h$-forms
on $U$ and denote the symmetric product of $S\in S^h(U)$ and $T\in
S^k(U)$ by $S\cdot T$. An element $L$ of $S^h(U)$ has a local
expression
\[
  L=L_{i_1\dots i_h}dx^{i_1}\dots dx^{i_h},
  \]
where the coefficients $L_{i_1\dots i_h}$ are smooth functions, which
are totally symmetric in the indices $i_1,\dots,i_h$. We
then define the $k$th order derivative of $L$ to be the
symmetric form of order $h+k$ given by
\[
 \delta^k(L) := \frac{\partial^k L_{i_1 \dots i_h}}{\partial
   x^{i_{h+1}}\dots\partial
    x^{i_{h+k}}}dx^{i_1}\dots dx^{i_h}dx^{i_{h+1}}\dots dx^{i_{h+k}}.
     \]
The definition depends on the choice of the local coordinates.

We have the following.

\begin{lemma}\label{l:contact}
Let $f : X \to \mathcal{Q}$ and $\hat{f} : X \to \mathcal{Q}$ be
smooth maps. Then, $f$ and $\hat{f}$ have analytic contact of order
$k$ at $p_0$ if and only if, for every local coordinate system $\{U;
x^1,x^2\}$ about $p_0$, there exist symmetric forms $\rho_r \in
S^r(U)$ such that
\begin{equation}\label{contact-condition}
  \delta^r(\hat{F})_{|p_0} = \sum_{h=0}^r \binom{r}{h} {\rho_{r-h}}_{|p_0}
   \delta^h({F})_{|p_0} \quad (r=0,\dots,k),
    \end{equation}
for arbitrary $F,\, \hat{F} : U \to \R^6$ such that
\[
 f(p) = [F(p)], \quad \hat{f}(p) = [\hat{F}(p)], \quad \text{for each}
  \quad p\in U.
   \]

\end{lemma}

\begin{proof}
Let $\{U; x^1,x^2\}$ be a coordinate system about $p_0$. As $G$ acts
transitively on $\mathcal{Q}$, we may assume that
\[
 f(p_0)=\hat{f}(p_0)=[e_0].
  \]
The map
\[
  y=(y^1,\dots,y^4)\in \R^4 \mapsto [X_0(y)] \in \mathcal{Q}
  \]
defined by
\begin{equation}
  X_0(y) = \Bigl(1,y,\displaystyle\frac{1}{2}(y,y)\Bigr)^T,
    \end{equation}
is a local coordinate system of $\mathcal{Q}$ centered at $[e_0]$.
Then, there exist an open neighborhood  $U'\subset U$ of $p_0$ and
smooth maps $h,\hat{h} : U'\to \R^4$ such that
\[
 f_{|U'}=[(X_0\circ h)],\quad
  \hat{f}_{|U'}=[(X_0\circ \hat{h})].
    \]
Thus, $f$ and $\hat{f}$ have analytic contact of order $k$ at $p_0$
if and only if the maps $\xi :=X_0\circ h$ and $\zeta :=X_0\circ \hat{h}$
satisfy
\begin{equation}\label{contact}
  \delta^r (\xi)_{|p_0}=\delta^r (\zeta)_{|p_0},\quad (r=0,\dots,k).
     \end{equation}
If $F$ and $\hat{F}$ are lifts of $f$ and $\hat{f}$, respectively, we
can write
\begin{equation}\label{F-xi}
  F = F^0 \xi,\quad \hat{F}= \hat{F}^0 \zeta
   \end{equation}
for smooth functions $F^0, \hat{F}^0 : U'\to \R$. Now, using \eqref{contact} and
\eqref{F-xi}, we compute
\begin{eqnarray*}
 \delta^r(\hat{F})_{|p_0} &=& \sum_{h=0}^r \binom{r}{h} \delta^{h}(\hat{F}^0)_{|p_0}
   \delta^{r-h}({\zeta})_{|p_0} =
\sum_{h=0}^r \binom{r}{h} \delta^{h}(\hat{F}^0)_{|p_0}
   \delta^{r-h}({\xi})_{|p_0} \\
&=&\sum_{h=0}^r \binom{r}{h} \delta^{h}(\hat{F}^0)_{|p_0}
   \delta^{r-h}\Bigl(\frac{F}{F^0}\Bigr)_{|p_0} \\
&=& \sum_{h=0}^r \sum_{m=0}^{r-h} \binom{r}{h} \binom{r-h}{m}
  \delta^{h}(\hat{F}^0)_{|p_0} \delta^{r-h-m}\Bigl(\frac{1}{F^0}\Bigr)_{|p_0}
\delta^{m}(F)_{|p_0}\\
&=& \sum_{m=0}^r \sum_{h=0}^{r-m} \binom{r}{m} \binom{r-m}{h}
  \delta^{h}(\hat{F}^0)_{|p_0} \delta^{r-m-h}\Bigl(\frac{1}{F^0}\Bigr)_{|p_0}
\delta^{m}(F)_{|p_0} \\
&=& \sum_{m=0}^r \binom{r}{m} \delta^{r-m}\Bigl(\frac{\hat{F}^0}{F^0}\Bigr)_{|p_0}
   \delta^{m}({F})_{|p_0},
  \end{eqnarray*}
and hence the conditions \eqref{contact-condition}.

Conversely, if conditions \eqref{contact-condition} hold true for
arbitrary lifts $F$ and $\hat{F}$, by choosing $F = \xi$ and $\hat{F} = \zeta$,
we get that $f$ and $\hat{f}$ have analytic contact of order $k$.
\end{proof}

\begin{cor}\label{c:contact}
In particular, we have proved that:

\begin{itemize}
\item $f$ and $\hat{f}$ have first order contact  if and only if
\begin{equation}\label{first-order}
 \hat{F}(p_0) = \rho_0(p_0) F(p_0), \quad
  \delta(\hat{F})_{|p_0} = \rho_1(p_0) F(p_0) + \rho_0(p_0)\delta(F)_{|p_0}.
   \end{equation}

\item $f$ and $\hat{f}$ have second order contact if and only if \eqref{first-order}
holds and
\begin{equation}\label{second-order}
  \delta^2(\hat{F})_{|p_0} = \rho_2(p_0) F(p_0) +
2\rho_1(p_0) \delta(F)_{|p_0} + \rho_0(p_0)\delta^2(F)_{|p_0}.
   \end{equation}

\item $f$ and $\hat{f}$ have third order contact if and only if
\eqref{first-order} and \eqref{second-order} hold and
\begin{eqnarray}\label{third-order}
  \delta^3(\hat{F})_{|p_0} & = & \rho_3(p_0) F(p_0) +
 3 \rho_2(p_0) \delta(F)_{|p_0} \\
  & &\mbox{}+  3\rho_1(p_0) \delta^2(F)_{|p_0} + \rho_0(p_0)\delta^3(F)_{|p_0}.\nonumber
   \end{eqnarray}
\end{itemize}
\end{cor}

\subsection{Conformal deformation}
In analogy with the characterization of conformal deformation of surfaces in
the conformal 3-sphere \cite{Mtrieste} (see also \cite{J}, \cite{JM}, \cite{MNtohoku}),
we can state the following result.

\begin{prop}\label{p:deformation}
Let $f, \hat{f} : X \to \Q$ be smooth spacelike immersions of the oriented
surface $X$ into the Lie quadric $\Q$, viewed as a homogeneous space of the
conformal group $G$.
Then, the following statements hold true:
\begin{enumerate}


\item The immersions $f$ and $\hat{f}$ are first order conformal deformation of
each other if and only if there exist first order frame fields $A$ and
$\hat{A}$ along $f$ and $\hat{f}$, respectively, such that
\begin{equation}\label{add-1def}
  {\hat{A}}^\ast {\omega}^1_0 = A^\ast{\omega}^1_0,\quad
  {\hat{A}}^\ast {\omega}^2_0 = A^\ast{\omega}^2_0,
   \end{equation}
where $\omega$ is the Maurer--Cartan form of $G$.

\item The immersions $f$ and $\hat{f}$ are second order conformal deformation of
each other if and only if there exist second order frame fields $A$ and
$\hat{A}$ along $f$ and $\hat{f}$, respectively, such that
\begin{equation}\label{add-2def}
\begin{array}{lll}
  {\hat{A}}^\ast {\omega}^1_0 = A^\ast{\omega}^1_0, &
  {\hat{A}}^\ast {\omega}^2_0 = A^\ast{\omega}^2_0, &
  {\hat{A}}^\ast {\omega}^2_1 = A^\ast{\omega}^2_1,\\
  {\hat{A}}^\ast {\omega}^0_0 = A^\ast{\omega}^0_0, &
  {\hat{A}}^\ast {\omega}^a_i = A^\ast{\omega}^a_i,& a=3,4;\,i= 1,2.
\end{array}
  \end{equation}

\item The immersions $f$ and $\hat{f}$ are third order conformal deformation of
each other if and only if there exist second order frame fields $A$ and
$\hat{A}$ along $f$ and $\hat{f}$, respectively, such that
 \[
  {\hat{A}}^\ast {\omega} = A^\ast{\omega}.
   \]
Thus any smooth immersion $f : X \to \mathcal{Q}$ is rigid to third order.

\end{enumerate}
\end{prop}

\begin{proof}[Sketch of the proof]

(1) Suppose that $f$ and $\hat{f}$ are first order conformal
deformations of each other. Then $D : U \to G$ exists so that
$\hat{f}$ and $D(p)f$ have first order analytic contact at $p$, for
each $p\in U$. Let $A : U \to G$ be a first order frame field along
$f$ and define $\hat{A} : U \to G$ by $\hat{A}(p) = D(p)A(p)$, for
each $p\in U$. Then $\hat{A}$ is a frame field along $\hat{f}$ and
${A}' = D(p_0)A : U \to G$ is a frame field along $D(p_0)f$, for
each $p_0 \in U$. According to \eqref{first-order} in Corollary
\ref{c:contact}, we have
\begin{eqnarray}
 {\hat{A}}_0 (p_0) & = & \rho_0(p_0)  {A'_0}(p_0)\label{1st-1}\\
  {d \hat{A}_0}_{|p_0} &=& \rho_0(p_0) {d A'_0}_{|p_0} +\rho_1(p_0)
   {A'_0}(p_0).\label{1st-2}
    \end{eqnarray}
Equation \eqref{1st-1} yields
\begin{equation}\label{ro-zero=1}
  \rho_0 =1,
 \end{equation}
since $\hat{A}$ and $A'$ agree at $p_0$. Now, the structure equations of $G$ imply
\begin{equation}\label{1-structure}
 d\hat{A}_0 = \hat{\theta}^J_0 \hat{A}_J, \quad
  dA'_0 = {\theta}^J_0 A'_J, \quad (J = 0,\dots,5),
  \end{equation}
where we have set $\theta = A^\ast\omega$ and $\hat\theta = \hat{A}^\ast\omega$.
Substituting \eqref{1-structure} into \eqref{1st-2} yields
\begin{equation}\label{eqns-1order}
  \rho_1 = (\hat{\theta}^0_0 -\theta^0_0), \quad
  \hat{\theta}^l_0 = \theta^l_0,\quad 1\leq l\leq 4.
    \end{equation}
Since $p_0$ has been chosen arbitrarily,
equations \eqref{eqns-1order} are identically satisfied on $U$.
Thus $\hat{A}$ is a
first order frame along $\hat{f}$ and
conditions \eqref{add-1def} are satisfied.

 Conversely, suppose \eqref{add-1def} hold for first order frame
fields $A$ and $\hat{A}$ along $f$ and $\hat{f}$, respectively. Then
define $D : U \to G$ by
\[
  D(p) = \hat{A}(p) A^{-1}(p), \quad \text{for each} \quad p\in U.
   \]
By \eqref{ro-zero=1}, \eqref{1-structure} and \eqref{eqns-1order}, we see that
\eqref{1st-1} and \eqref{1st-2} hold. So, $D$ induces a first order conformal
deformation.

\vskip0.1cm

(2) Retain the notations of (1) and suppose that $f$ and $\hat{f}$
are second order conformal deformations of each other. Then
$\hat{f}$ and $D(p)f$ have second order analytic contact at $p$, for
each $p\in U$. Let $A : U \to G$ be a second order frame field along
$f$ and $\hat{A}$ and $A'$ be as in (1). We have to show that
$\hat{A}$ defines a second order frame field along $\hat{f}$ such
that
\begin{equation}\label{2order-conds}
\begin{array}{lll}
  \hat{\theta}^1_0 = \theta^1_0, &
  \hat{\theta}^2_0 = \theta^2_0, &
  \hat{\theta}^2_1 = \theta^2_1,\\
  \hat{\theta}^0_0 = \theta^0_0, &
  \hat{\theta}^a_i = \theta^a_i &
 (a= 3,4; i= 1,2).
\end{array}
  \end{equation}
By Corollary \ref{c:contact} and the discussion in part (1), we know that
the frame fields $\hat{A}$ and $A'$ must satisfy
\eqref{1st-1}, \eqref{1st-2} and
 \begin{equation}\label{2nd-3}
  \delta^2(\hat{A}_0)_{|p_0} = \rho_2(p_0) A'_0(p_0) +
    2(\hat{\theta}^0_0 -\theta^0_0)_{|p_0} d {A'_0}_{|p_0} +
    \delta^2(A'_0)_{|p_0}.
      \end{equation}
Writing out \eqref{2nd-3}, using the structure equations
\eqref{1-structure}, the equations $\hat{\theta}^i_0 ={\theta}^i_0$,
$i=1,2$, and $\hat{\theta}^a_0 ={\theta}^a_0 =0$, $a=3,4$, and the
fact that $\hat{A}_0(p_0) = A'_0(p_0)$, we find
\begin{eqnarray}
 \rho_2 &=& \delta(\hat{\theta}^0_0 - \theta^0_0) +
  (\hat{\theta}^0_0 - \theta^0_0)^2 + \theta^1_0(\hat{\theta}^0_1 -{\theta}^0_1)
  + \theta^2_0(\hat{\theta}^0_2 -{\theta}^0_2),\\
  0 &=& \theta^2_0 (\hat{\theta}^1_2 -{\theta}^1_2)
  -\theta^1_0 (\hat{\theta}^0_0 - \theta^0_0), \\
   0 &=& \theta^1_0 (\hat{\theta}^2_1 -{\theta}^2_1)
  -\theta^2_0 (\hat{\theta}^0_0 - \theta^0_0), \\
    0 &=& \theta^1_0 (\hat{\theta}^a_1 - \theta^a_1) +
  \theta^2_0 (\hat{\theta}^a_2 - \theta^a_2) \quad (a = 3,4).
    \end{eqnarray}
From these equations it follows that
\[
 \hat{\theta}^2_1 = {\theta}^2_1, \quad \hat{\theta}^0_0 - \theta^0_0, \quad
  \hat{\theta}^a_i = \theta^a_i \quad (i=1,2; a=3,4).
   \]
Thus $\hat{A}$ is a second
order frame along $\hat{f}$ and the conditions \eqref{2order-conds} are satisfied.

 Conversely, suppose \eqref{2order-conds} hold for second order frame
fields $A$ and $\hat{A}$ along $f$ and $\hat{f}$, respectively. As
above, define $D : U \to G$ by
\[
  D(p) = \hat{A}(p) A^{-1}(p), \quad \text{for each} \quad p\in U.
   \]
By reversing the arguments above. we see that \eqref{1st-1}, \eqref{1st-2}
and \eqref{2nd-3} are satisfied. So, $D$ induces a second order conformal
deformation of $f$ and $\hat{f}$.

\vskip0.1cm
As for (3), writing out \eqref{third-order}, one can prove after some lengthy
computations that $\theta = A^{-1}dA = \hat{A}^{-1}d\hat{A} = \hat{\theta}$.
By the Cartan--Darboux rigidity theorem, one then have $dD_{|p} = 0$,
for every $p \in U$.

\end{proof}

\begin{ex}[Isothermic surfaces]\label{e:isothermic}

Let $X$ be an oriented 2-di\-mensio\-nal manifold and let
$f : X \to \Q \subset \RP^5$
be a smooth spacelike conformal immersion. On $X$, consider the
unique complex structure defined by the given orientation and the
conformal structure induced by $f$. We recall the following.

\begin{defn}
The immersion $f: X \to \Q$ is \textit{isothermic} if there exists a
non-zero holomorphic quadratic differential $Q$ on $X$ and a section
$S$ of $\mathcal{N}(X)\otimes K_X$ such that $\mathcal{A}_{|q} =
Q_{|q}\otimes S_{|q}$, for each $q\in X$ such that $Q_{|q} \neq
0$.\footnote{Near any point $q\in X$ such that $Q_{|q} \neq 0$,
there is a complex parameter $z = x +iy$ such that $Q = dz\,dz$.
Further, if $z$ and $\tilde{z}$ are two such coordinates, then $dz =
\pm d\tilde{z}$. If $f$ is isothermic, then $(x,y)$ are principal,
isothermal local coordinates, and the immersion $f$ is isothermic in
the classical sense. The converse holds only locally. If $f$ is not
totally umbilical, then the holomorphic differential $Q$ is uniquely
defined up to a non-zero constant factor (see \cite{Pa1988},
\cite{Ber2001} and \cite{HJlibro}, Chapter 5).} Locally, this means
that
\[
 \mathcal{A} = r^a Q A_a \otimes A_0
  \]
for real-valued smooth functions $r^a$, $a = 3,4$.

\end{defn}

We now discuss the nonlinear partial differential equation governing
iso\-thermic surfaces; we mainly follow \cite{BHPP} (see also
\cite{HJlibro} and \cite{Burs}). First, note that the covariant
derivative $D^\nu$ on the conformal normal bundle of an isothermic
immersion is flat. Further, fix coordinates $z = x + iy$ such that
$Q = dz\,dz$ and choose a flat trivialization $(A_3, A_4)$ of the
conformal normal bundle. Then there exists a unique cross section $A
: X \to P_2(f)$ such that
\begin{equation}\label{connection-form}
 A^{-1}dA=\left(
 \begin{array}{cccccc}
  0&\chi_1&\chi_2&\tau_1&\tau_2&0\\
  dx&0&0&-k_1dx&k_2dx&\chi_1\\
   dy&0&0&k_1dx&-k_2dx&\chi_2\\
    0&k_1dx&-k_1dx&0&0&\tau_1\\
     0&k_2dx&-k_2dy&0&0&\tau_2\\
       0&dx&dy&0&0&0
         \end{array}\right),
          \end{equation}
where $k_1,k_2 : X \to \R$ are smooth functions and $\chi_1,
\chi_2$, $\tau_1,\tau_2$ 1-forms to be determined. (In general, $A$
depends on the polarization $Q$, since totally umbilical immersions
are allowed. This dependence disappears if $f$ is not totally
umbilical). From the Maurer--Cartan equations, we have
\begin{equation}\label{chi-tau}
 \begin{aligned}
  \chi_1 & = \frac{1}{2}(u - \|k\|)dx + \psi dy, &
     \chi_2 & = \psi dx -\frac{1}{2}(u + \|k\|)dy \\
      \tau_1 & = (k_1)_xdx - (k_1)_ydy, &
     \tau_2 & = -(k_2)_xdx + (k_2)_ydy, \\
    \end{aligned}
    \end{equation}
where
$$
  \|k\| := (k_1)^2 -(k_2)^2
   $$
and $\psi$, $u$ are smooth functions satisfying
\[
 \begin{aligned}
  & (k_1)_{xy}  =  \psi k_1, &
     & (k_2)_{xy}  = \psi k_2, \\
   & u_x   = -2 ( \psi_y + \|k\|_x), &
     & u_y  = 2 ( \psi_x + \|k\|_y). \\
     \end{aligned}
       \]
The above equations are compatible if and only if $k_1$, $k_2$ and
$\psi$ are solutions of the \textit{vector Calapso equation}
\[
 (k_1)_{xy} =\psi k_1,\quad  (k_2)_{xy}  = \psi k_2,
  \quad \Delta \psi = -2 \|k\|_{xy}.
    \]

Conversely, if we start with a solution $(k_1,k_2,\psi)$ of the
vector Calapso equation, then the differential 1-form
\[
 \upsilon =-2 \left( \psi_y  +  \|k\|_x \right)dx +2 \left(
   \psi_x  +  \|k\|_y \right)dy
    \]
is closed. Let $u : X \to \R$ be a primitive of $\upsilon$, and
define the 1-forms $\chi_1,\chi_2$ and $\tau_1, \tau_2$ as in
\eqref{chi-tau}. Next, define a $\mathfrak{g}$-valued form $\alpha$
as in \eqref{connection-form}; $\alpha$ satisfies the Maurer--Cartan
equation and integrates (locally) to a smooth map $A : X \to G$ such
that $\alpha = A^{-1}dA$. Then $f = [A_0] : X \to \Q$ is an
isothermic immersion. Since $u$ is defined up to a constant, for
each $\lambda \in \R$, there exists a frame $A_\lambda$ and a
corresponding isothermic immersion $f_\lambda$. Define
\[
 D_\lambda = A_\lambda A^{-1}  : X \to G.
  \]
It is easy to see that $f$ and $f_\lambda$ are second order
deformation of each other with respect to the infinitesimal
displacement $D_\lambda$. The isothermic immersion $f_\lambda$
coincides with the classical T-transform (spectral deformation) of
the isothermic immersion $f$, introduced independently by L. Bianchi
and P. Calapso at the turn of the 20th century. The displacement
$D_\lambda$ induces the T-transformation of isothermic surfaces.

It follows that isothermic surfaces are (locally) deformable of
second order and allow 1-parameter families of second order
deformations. Such families correspond to the solutions of the
vector Calapso equation which in turn is equivalent to the
Gauss-Codazzi equations and arises as integrability condition of a
linear differential system containing a free parameter. For more on
the relations with the theory of integrable systems and the
classical theory of transformations of isothermic surfaces,
including the T-transformation, we refer the reader to
\cite{HJlibro}, Chapter 5, and the references therein.

\end{ex}

Consider the Pfaffian differential systems $\mathfrak{I}_1$,
$\mathfrak{I}_2$ and $\mathfrak{I}_3$ with independence condition
$\alpha^1_0 \wedge \alpha^2_0 \neq 0$, introduced in Section
\ref{ss:conf-space}. According to Proposition \ref{p:deformation},
Example \ref{e:isothermic}, and the preceding discussion, we can
state the following.

\begin{thm}

{\rm (1)} The integral manifolds of the Pfaffian differential system
$(\mathfrak{I}_1, \alpha^1_0 \wedge \alpha^2_0)$
arise as maps
\[
   \mathfrak{d} : X \to \D, \, q \mapsto (\tau_f(q), \tau_{\hat{f}}(q), D(q)),
    \]
where $f, \hat{f} : X \to \mathcal{Q}$ are spacelike immersions which are first
order deformations of each other and $D : X \to G$ is the
infinitesimal displacement of the deformation.

{\rm (2)} The integral manifolds of the Pfaffian differential system
$(\mathfrak{I}_2, \alpha^1_0 \wedge \alpha^2_0)$
arise as maps
\[
   \mathfrak{d} : X \to \D, \, q \mapsto (\tau_f(q), \tau_{\hat{f}}(q), D(q)),
    \]
where $f, \hat{f} : X \to \mathcal{Q}$ are spacelike immersions which are second
order deformations of each other and
$D : X \to G$ is the infinitesimal displacement of the deformation.

{\rm (3)} The integral manifolds of the Pfaffian differential system
$(\mathfrak{I}_3, \alpha^1_0 \wedge \alpha^2_0)$
arise as maps
\[
   \mathfrak{d} : X \to \D, \, q \mapsto (\tau_f(q), \tau_{\hat{f}}(q), D(q)),
    \]
where $f, \hat{f} : X \to \mathcal{Q}$ are spacelike isothermic immersions which
are T-transforms of each other and $D : X \to G$ defines the T-transformation.
\end{thm}

It follows that the study of second order conformal deformations of
spacelike surfaces reduces to the study of the integral manifolds of
$(\mathfrak{I}_2, \alpha^1_0 \wedge \alpha^2_0)$. In particular, the
study of the class of deformable surfaces given by isothermic
surfaces reduces to the study of the integral manifolds of
$(\mathfrak{I}_3,\alpha^1_0 \wedge \alpha^2_0)$.

\section{The exterior differential system of a deformation}\label{s:eds}

In this section, we undertake the study of the PDS of second order conformal deformations.

Let $(\alpha, \beta)$ be the $\mathfrak{g}\times\mathfrak{g}$-valued 1-form
on ${\mathcal{D}}$ obtained by pulling-back
the Maurer-Cartan form $(\omega, \Omega)$ of $G\times G$ with respect to a
local section of the projection $\pi_{\mathcal{D}} : G\times G \to {\mathcal{D}}$.
Let $\alpha^1=\alpha^1_0$,
$\alpha^2=\alpha^2_0$ and set
\begin{equation}\label{Generators1}
 \left\{
\begin{array}{lll}
\eta^1 =\alpha^0_0 - \beta^0_0, &\eta^2 = \alpha^1_0-\beta^1_0, &\eta^3 =\alpha^2_0-\beta^2_0,\\
\eta^4 =\alpha^3_0, & \eta^5 =\alpha^4_0, & \eta^6 = \beta^3_0, \\
\eta^7 =\beta^4_0,& \eta^8 =\alpha^2_1 -\beta^2_1, &\eta^9  =\alpha^3_1-\beta^3_1,\\
\eta^{10} =\alpha^3_2 -\beta^3_2,& \eta^{11} =\alpha^4_1 -\beta^4_1, &\eta^{12}  =\alpha^4_2-\beta^4_2.
     \end{array}\right.
      \end{equation}

The Pfaffian differential system $(\mathfrak{I}_2,
\alpha^1\wedge\alpha^2)$ differentially generated by the 1-forms
$\eta^1, \dots,\eta^{12}$ with independent condition
\[
 \alpha^1 \wedge \alpha^2 \neq 0
  \]
is called the \textit{differential system of a deformation.}

\begin{remark}
The definition of $\mathfrak{I}_2$ is independent of the local
sections of $\pi_{\mathcal{D}}$. The integral manifolds of
$(\mathfrak{I}_2,\alpha^1\wedge\alpha^2)$ are the two-dimensional
immersed submanifolds
\[
 \mathfrak{d} = \left([A_0\wedge A_1 \wedge A_2], [B_0\wedge B_1
   \wedge B_2], BA^{-1}\right) : X \to \mathcal{D},
     \]
where $(A,B) : X \to G\times G$ is an integral manifold of the
Pfaffian differential system $\pi_{\mathcal{D}}^\ast
(\mathfrak{I}_2)$ defined on $G\times G$ and $f=[A_0]$, $\hat{f} =
[B_0]$ are second order deformations of each other with respect to
$D = BA^{-1}$, i.e., $\hat{f}$ and $D(q)\cdot f$ have second order
analytic contact at $q$, for each $q\in X$.
\end{remark}

\vskip0.3cm

Using the Maurer-Cartan equations, we compute, modulo the algebraic
ideal generated by $\eta^1,\dots,\eta^{12}$, the \textit{quadratic
equations} of the system
\begin{equation} \label{QE1}
\left\{ \begin{array}{lcl}
  d\eta^1&\equiv& - (\alpha^0_1 -\beta^0_1)\wedge\alpha^1 - (\alpha^0_2 -\beta^0_2)\wedge\alpha^2,\\
    d\eta^2&\equiv& d\eta^3 \equiv 0,\\
     d\eta^4&\equiv& d\eta^6 \equiv - \alpha^3_1 \wedge\alpha^1 - \alpha^3_2 \wedge\alpha^2,\\
    d\eta^5 &\equiv& d\eta^7\equiv - \alpha^4_1 \wedge\alpha^1 - \alpha^4_2 \wedge\alpha^2,\\
\end{array}\right.
\end{equation}
 \begin{equation} \label{QE2}
\left\{ \begin{array}{lcl}
  d\eta^8&\equiv&  (\alpha^0_1 -\beta^0_1)\wedge\alpha^2 - (\alpha^0_2 -\beta^0_2)\wedge\alpha^1,\\
  d\eta^9&\equiv&  -(\alpha^0_3 -\beta^0_3)\wedge\alpha^1 - (\alpha^4_3 -\beta^4_3)\wedge\alpha^4_1,\\
    d\eta^{10}&\equiv& -(\alpha^0_3 -\beta^0_3)\wedge\alpha^2 - (\alpha^4_3 -\beta^4_3)\wedge\alpha^4_2,\\
    d\eta^{11}&\equiv& (\alpha^0_4 -\beta^0_4)\wedge\alpha^1 - (\alpha^4_3 -\beta^4_3)\wedge\alpha^3_1,\\
     d\eta^{12}&\equiv&  (\alpha^0_4 -\beta^0_4)\wedge\alpha^2 - (\alpha^4_3 -\beta^4_3)\wedge\alpha^3_2.
\end{array}\right.
\end{equation}

From this, we see that the differential ideal $\mathfrak{I}_2$ is
algebraically generated by $\eta^1,\dots,\eta^{12}$ and the 2-forms
\begin{equation}\label{QuadraticGenerators}
  \left\{
\begin{array}{lcl}
\Omega^1 &=& - (\alpha^0_1 -\beta^0_1)\wedge\alpha^1 - (\alpha^0_2 -\beta^0_2)\wedge\alpha^2,\\
 \Omega^2 &=& - \alpha^3_1 \wedge\alpha^1 - \alpha^3_2 \wedge\alpha^2,\\
  \Omega^3 &=& - \alpha^4_1 \wedge\alpha^1 - \alpha^4_2 \wedge\alpha^2,\\
   \Omega^4 &=& (\alpha^0_1 -\beta^0_1)\wedge\alpha^2 - (\alpha^0_2 -\beta^0_2)\wedge\alpha^1,\\
    \Omega^5 &=& -(\alpha^0_3 -\beta^0_3)\wedge\alpha^1 - (\alpha^4_3 -\beta^4_3)\wedge\alpha^4_1,\\
     \Omega^6 &=& -(\alpha^0_3 -\beta^0_3)\wedge\alpha^2 - (\alpha^4_3 -\beta^4_3)\wedge\alpha^4_2,\\
     \Omega^7 &=& (\alpha^0_4 -\beta^0_4)\wedge\alpha^1 - (\alpha^4_3 -\beta^4_3)\wedge\alpha^3_1,\\
     \Omega^8 &=&  (\alpha^0_4 -\beta^0_4)\wedge\alpha^2 - (\alpha^4_3 -\beta^4_3)\wedge\alpha^3_2.
      \end{array}\right.
       \end{equation}

\begin{remark}
The Pfaffian differential system
$(\mathfrak{I}_2,\alpha^1\wedge\alpha^2)$ is \textit{linear}, in the
sense that
\[
  d\left\{\eta^1,\dots,\eta^{12} \right\} \subset
   \left\{\eta^1,\dots,\eta^{12}, \alpha^1,\alpha^2 \right\},
    \]
where $\left\{\eta^1,\dots,\eta^{12} \right\}$ and
$\left\{\eta^1,\dots,\eta^{12}, \alpha^1,\alpha^2 \right\}$ denote,
respectively, the ideals generated algebraically by
$\eta^1,\dots,\eta^{12}$ and $\eta^1,\dots,\eta^{12},
\alpha^1,\alpha^2$. If a Pfaffian differential system is linear,
then its integral elements at a point are determined by
inhomogeneous linear equations. In the literature, linear systems
are also called quasi-linear systems, systems in good form, or
systems in normal form.

\end{remark}

We now discuss the involutiveness of the system. For this, we
consider the basis
\[
   (\alpha^1,\,\alpha^2,\,\eta^j,\, \alpha^a_i,\,\alpha^0_I-\beta^0_I,\, \alpha^4_3-\beta^4_3 )
     \]
($j=1.\dots,12$; $i=1,2$; $a=3,4$; $I = 1,2,3,4$) for the 1-forms on $\mathcal{D}$ and denote
by
\[
 \left(\frac{\partial}{\partial \alpha^1}, \frac{\partial}{\partial \alpha^2},\frac{\partial}{\partial
   \eta^j},\frac{\partial}{\partial {\alpha^a_i}},
      \frac{\partial}{\partial
      {(\alpha^0_I-\beta^0_I)}}, \frac{\partial}{\partial {(\alpha^4_3-\beta^4_3)}}\right)
      \]
its dual basis.

A 1-dimensional integral element of the system is of the form $E_1 = [V]$, where
\begin{eqnarray*}
  V &=& a_{i}\frac{\partial}{\partial
   \alpha^i}+b_1\frac{\partial}{\partial {\alpha^3_1}}+b_2\frac{\partial}{\partial {\alpha^3_2}}
+b_3\frac{\partial}{\partial {\alpha^4_1}}+b_{4}\frac{\partial}{\partial {\alpha^4_2}} \\
& &\mbox{}+b_{4+I}\frac{\partial}{\partial
      {(\alpha^0_I-\beta^0_I)}}
     +b_9\frac{\partial}{\partial {(\alpha^4_3-\beta^4_3)}}
    \end{eqnarray*}
is a general vector in the space $\eta^j =0$, $j=1,\dots,12$.
Thus, the manifold of 1-dimensional integral elements
$\mathcal{V}_1\cong \mathcal{D}\times \RP^{10}$.
Moreover, $E_1$ is admissible if and only if $(a_1)^2+(a_2)^2\neq 0$.

The \textit{polar equations} of a given $E_1\in\mathcal{V}_1$ are
$\eta^j=0$ ($j =1,\dots,12$) and
\[
i_{V}\Omega^\beta =0 \qquad (\beta =1,\dots,8),
   \]
which read
\begin{equation}\label{polareqs}
\begin{aligned}
  -b_5\alpha^1 + a_1(\alpha^0_1 -\beta^0_1) -b_6\alpha^2 + a_2(\alpha^0_2 -\beta^0_2) &=0 ,\\
 -b_1\alpha^1 +a_1\alpha^3_1 - b_2\alpha^2 + a_2 \alpha^3_2 &=0,\\
-b_3\alpha^1 + a_1\alpha^4_1 - b_4\alpha^2 + a_2 \alpha^4_2 &=0,\\
b_5\alpha^2 - a_2(\alpha^0_1 -\beta^0_1) -b_6\alpha^1 + a_1(\alpha^0_2 -\beta^0_2) &=0 ,\\
-b_7\alpha^1 + a_1(\alpha^0_3 -\beta^0_3) -b_9\alpha^4_1 + b_3(\alpha^4_3 -\beta^4_3) &=0 ,\\
-b_7\alpha^2 + a_2(\alpha^0_3 -\beta^0_3) -b_9\alpha^4_2 + b_4(\alpha^4_3 -\beta^4_3) &=0 ,\\
b_8\alpha^1 - a_1(\alpha^0_4 -\beta^0_4) -b_9\alpha^3_1 + b_1(\alpha^4_3 -\beta^4_3) &=0 ,\\
b_8\alpha^2 - a_2(\alpha^0_4 -\beta^0_4) -b_9\alpha^3_2 + b_2(\alpha^4_3 -\beta^4_3) &=0.\\
\end{aligned}
\end{equation}

If
\[
 b_9 \left[(a_1)^2-(a_2)^2\right] \left[a_1(b_2b_7 +b_4b_8) -a_2(b_1b_7 +b_3b_8) +
   b_9(b_2b_3 -b_1b_4)\right] \neq 0,
    \]
the polar equations are linearly independent and the polar space
$H(E_1)$ is 3-dimensional. Thus $c_0 = \dim \mathfrak{I}_2^{(1)} =
12$ and $c_1 = \codim H(E_1) = 20$. Further, the variety of
2-dimensional integral elements over a point $x\in \mathcal{D}$ has
dimension $10$. Hence $c_0 + c_1 = \dim G_2(T_x\mathcal{D}) - \dim
\mathcal{V}_2(x) = 32$, and Cartan's test applies. The Cartan
characters $s_k = c_k -c_{k-1}$, $k = 0,1,2$, are then computed to
be $s_0 = 12$, $s_1 = 8$, $s_2 = 1$.

\vskip0.1cm

Summarizing, we can state the following.

\begin{prop}
The Pfaffian differential system
$(\mathfrak{I}_2,\alpha^1\wedge\alpha^2)$ is in involution and its
general solutions depend on one arbitrary function in two
variables.\footnote{Considering that a generic surface in a
($2+r$)-dimensional space may be locally given as graph of $r$
arbitrary functions of two variables, deformable surfaces in
compactified Minkowski 4-space are then exceptional.} The singular
solutions of the system correspond to the points of the reducible
variety defined by the equation
\[
  b_9 \left[(a_1)^2-(a_2)^2\right] \left[a_1(b_2b_7 +b_4b_8) -a_2(b_1b_7 +b_3b_8) +
   b_9(b_2b_3 -b_1b_4)\right] = 0.
    \]

In particular, isothermic surfaces correspond to the points of the variety defined by
\[
 b_9  \left[(b_7)^2 + (b_8)^2 + (b_9)^2\right] = 0.
\]

\end{prop}

\begin{remark}
Note that isothermic surfaces satisfy the additional equations
\[
\alpha^4_3 -\beta^4_3 =0,\quad \alpha^0_3 -\beta^0_3 =0,\quad
\alpha^0_4 -\beta^0_4 =0,
      \]
and are then integral manifolds of the Pfaffian differential system
$(\mathfrak{I}_3,\alpha^1\wedge\alpha^2)$, differentially generated
by $I_3 \subset \Gamma (T^\ast\mathcal{D})$ (see Section
\ref{ss:conf-space}). A direct computation shows that the system
$(\mathfrak{I}_3,\alpha^1\wedge\alpha^2)$ is in involution and its
general solution depends on six arbitrary functions in one variable.

\end{remark}

\bibliographystyle{amsalpha}

\end{document}